\documentclass[11pt,twoside]{article}

\usepackage{a4wide}
\usepackage{amsfonts}
\usepackage{amssymb}
\usepackage{amsmath}
\usepackage{graphicx}
\usepackage{nopageno}
\usepackage{xcolor}

\newcommand{\ignore}[1]{}

\def\@begintheorem#1#2{\par\bgroup{\sc #1\ #2. }\it\ignorespaces}
\def\@opargbegintheorem#1#2#3{\par\bgroup{\sc #1\ #2\ (#3). } \it\ignorespaces}
\def\@endtheorem{\egroup}
\newtheorem{theorem}{Theorem}[section]
\newtheorem{corollary}[theorem]{Corollary}
\newtheorem{lemma}[theorem]{Lemma}

\newtheorem{example}[theorem]{Example}
\newtheorem{proposition}[theorem]{Proposition}
\newtheorem{definition}[theorem]{Definition}
\newcommand{\bt}[1]{\begin{theorem}\label{#1}}
\newcommand{\bc}[1]{\begin{corollary}\label{#1}}
\newcommand{\bl}[1]{\begin{lemma}\label{#1}}
\newcommand{\be}[1]{\begin{example}\label{#1}}
\newcommand{\bp}[1]{\begin{proposition}\label{#1}}
\newcommand{\ba}[1]{\begin{algorithm}\rm\label{#1}}
\newcommand{\bd}[1]{\begin{definition}\rm\label{#1}}{\normalsize }

\newcommand{\et}{\end{theorem}}
\newcommand{\ec}{\end{corollary}}
\newcommand{\el}{\end{lemma}}
\newcommand{\ee}{\end{example}}
\newcommand{\ep}{\end{proposition}}
\newcommand{\ed}{\end{definition}}
\newcommand{\epr}{{\ \vbox{\hrule\hbox{%
\vrule height1.3ex\hskip0.8ex\vrule}\hrule}}\\\par}

\def\Z{\mathbb{Z}}
\def\H{{\cal H}}
\def\M{{\cal M}}

\begin{document}

\title{\bf Matching Orderable and Separable Hypergraphs}

\author{
Shmuel Onn
\thanks{\small Technion - Israel Institute of Technology. Email: onn@technion.ac.il.}
\thanks{\small My manuscript has no associated data.}
}

\date{}

\maketitle

\begin{abstract}
A perfect matching in a hypergraph is a set of edges that partition the set of vertices.
We study the complexity of deciding the existence of a perfect matching in orderable and separable
hypergraphs. We show that the class of orderable hypergraphs is strictly contained in the class
of separable hypergraphs.
Accordingly, we show that for each fixed $k$, deciding perfect matching
for orderable $k$-hypergraphs is polynomial time doable, but for each fixed $k\geq 3$, it
is NP-complete for separable hypergraphs.

\end{abstract}

\section{Introduction}

A {\em $k$-hypergraph} on a finite set $V$ is a set $\H$ of $k$-subsets of $V$. The elements
$v\in V$ are called {\em vertices} and the sets $E\in\H$ are called {\em edges}.
A {\em perfect matching} in $\H$ is a subset $\{E_1,\dots,E_m\}\subseteq \H$ such that
$V=\uplus_{i=1}^m E_i$, that is, $V$ is the disjoint union of $E_1,\dots,E_m$.

The special case of perfect matchings in graphs (that is, $2$-hypergraphs) has numerous
applications and is one of the most studied combinatorial optimization problems \cite{Sch}.
The natural extension to $k$-hypergraphs for $k\geq 3$ has also been studied by several authors
and has  a variety of applications, such as for the Santa Claus allocation problem,
see e.g. \cite{AFS,Cyg,Hax} and the references therein.

In this article we consider the complexity of the problem
of deciding if a perfect matching exists. We assume throughout that $k$ is fixed.
Clearly a necessary condition for a perfect matching to exist is that $n=km$ is a
multiple of $k$ where $n:=|V|$.
For $k=1$ the problem is trivial as there is a perfect matching if and only if
$\H=\{\{v\}:v\in V\}$. For $k=2$, that is, graphs, it is well known
that the problem can be solved in polynomial time, see e.g. \cite{Sch}.

For $k=3$ the perfect matching problem, and very special cases of it including the
following two, are NP-complete, see \cite{GJ}. First, the {\em $3$-dimensional matching} problem,
where we are given a partition $V=V_1\uplus V_2\uplus V_3$ of $V$ into three $m$-sets,
the given hypergraph satisfies
$\H\subseteq\{E\subseteq V:|E\cap V_1|=|E\cap V_2|=|E\cap V_3|=1\}$,
and we need to decide if $\H$ has a perfect matching.

Second, the {\em $3$-partition} problem, where we are given a labeling $a:V\rightarrow\Z$ of vertices
by integers, which can even be encoded in unary,
and denoting $a(U):=\sum\{a(u):u\in U\}$ for any subset $U\subseteq V$,
the hypergraph is given by $\H=\{E\subseteq V:|E|=3,\ a(E)=0\}$,
and we need to decide if $\H$ has a perfect matching.

This hardness of the problem in general motivates various lines of investigation.
First, it is interesting to study algorithms for approximating a perfect matching. This line
is taken in \cite{Cyg} and the references therein, where effective approximation algorithms
for $3$-dimensional matching are given. Second, it is interesting to classify the
complexity of deciding perfect matching in subclasses of hypergraphs. This line is taken,
for example, in \cite{Hax}, where a criterion for perfect matching in {\em bipartite}
hypergraphs is derived.

Here we continue this second line of investigation and consider two classes
of hypergraphs defined below. These two classes arise as natural extensions to hypergraphs
of the well studied class of {\em threshold graphs}, which in itself has many applications
and interesting subclasses, see \cite{MP}.

First, a $k$-hypergraph is {\em orderable} if there is an {\em elimination order} $v_1,\dots,v_n$
of $V$, that is, an order where each vertex $v_i$ is either {\em dominating},
meaning that $E\in\H$ for every $k$-set $E$ with $v_i\in E\subseteq\{v_1,\dots,v_i\}$,
or {\em isolating}, meaning that $E\notin\H$ for every $k$-set $E$ with
$v_i\in E\subseteq\{v_1,\dots,v_i\}$. We prove the following theorem.

\bt{orderable}
For any fixed $k$, the following two statements hold:
\begin{enumerate}
\item
It can be decided in polynomial time if a given $k$-hypergraph is orderable;
\item
It can be decided in polynomial time if an orderable $k$-hypergraph has a
perfect matching.
\end{enumerate}
\et

Second, reminiscent of the $3$-partition problem above, a $k$-hypergraph is {\em separable}
if there is a labeling $a:V\rightarrow\Z$ of vertices by integers such that
$\H=\{E\subseteq V:|E|=k,\ a(E)\geq 0\}$. For $k=2$ this coincides with the well studied
class of threshold graphs mentioned above. We show in Proposition \ref{orderable_separable}
that the classes of orderable and separable hypergraphs coincide for $k=1,2$, but the former
is strictly contained in the latter for all $k\geq 3$. Thus, deciding the existence of a
perfect matching for this broader class is
expected to be harder, and we confirm it by proving the following theorem.

\bt{separable}
For any fixed $k$, the following two statements hold:
\begin{enumerate}
\item
It can be decided in polynomial time if a given $k$-hypergraph is separable;
\item
Deciding if a separable $k$-hypergraph has a perfect matching is
polynomial time doable for $k=1,2$, but is NP-complete for each $k\geq 3$, even
if $a:V\rightarrow\Z$ is encoded in unary.
\end{enumerate}
\et

As a rather general line of investigation, it would be interesting to identify meaningful classes
of hypergraphs which lie in between the classes of orderable and separable hypergraphs,
and to understand the complexity of deciding the existence of a perfect matching for such classes.

\section{Orderable hypergraphs}

We now prove our theorem about orderable hypergraphs.

\vskip.2cm\noindent{\bf Theorem \ref{orderable}}
{\em For any fixed $k$, the following two statements hold:
\begin{enumerate}
\item
It can be decided in polynomial time if a given $k$-hypergraph is orderable;
\item
It can be decided in polynomial time if an orderable $k$-hypergraph has a
perfect matching.
\end{enumerate}
}

{\em Proof of part 1.}
First, we claim that if $\H$ is orderable, and $v\in V$ satisfies that either $E\in\H$ for each
$k$-set $E$ with $v\in E\subseteq V$, or $E\notin\H$ for each $k$-set $E$ with
$v\in E\subseteq V$, then $\H$ has an elimination order where $v$ comes last. To see this,
suppose $v_1,\dots,v_n$ is any elimination order for $\H$ with $v=v_i$ and consider the
order $v_1,\dots,v_{i-1},v_{i+1},\dots,v_n,v_i$. We claim this is also an elimination order.
For $v_i=v$ the condition holds by assumption. Since $v_1,\dots,v_n$ is an elimination
order, for each other vertex $v_j$ we have that either $E\in\H$ for each $k$-set
$E$ with $v_j\in E\subseteq\{v_1,\dots,v_j\}$, or $E\notin\H$ for each $k$-set $E$ with
$v_j\in E\subseteq\{v_1,\dots,v_j\}$, and hence also either $E\in\H$ for each $k$-set
$E$ with $v_j\in E\subseteq\{v_1,\dots,v_j\}\setminus\{v_i\}$, or $E\notin\H$ for each
$k$-set $E$ with $v_j\in E\subseteq\{v_1,\dots,v_j\}\setminus\{v_i\}$.
So this is indeed an elimination order where $v$ comes last.

This implies that $\H$ is orderable if and only if there exists a vertex $v\in V$
for which the following hold: first, either $v$ is dominating, that is, $E\in\H$
for each $k$-set $E$ with $v\in E\subseteq V$, or $v$ is isolating, that is,
$E\notin\H$ for each $k$-set $E$ with $v\in E\subseteq V$; and second, the hypergraph
$\H':=\{E\in\H:E\subseteq V'\}$ on $V':=V\setminus\{v\}$ is also orderable.

So the algorithm proceeds recursively as follows. If $n:=|V|\leq k$ then $\H$ is orderable
by definition. If $n>k$ then we search for $v\in V$ such that either $E\in\H$ for each $k$-set
$E$ with $v\in E\subseteq V$, or $E\notin\H$ for each $k$-set $E$ with $v\in E\subseteq V$.
If there is no such $v$ then $\H$ is not orderable. If there is such $v$ then we define
$V'$ and $\H'$ as above and apply the algorithm recursively to $\H'$.
The running time of the algorithm is dominated by the number $t(n)$ of $k$-sets $E$ for
which we test if $E\in\H$. We show by induction on $n$ that $t(n)\leq n^{k+1}$.
For $n\leq k$ any $\H$ is orderable so $t(n)=0$. Suppose $n>k$. For each $v\in V$,
checking if $E\in\H$ for each $k$-set $E$ with $v\in E\subseteq V$, or $E\notin\H$ for each
$k$-set $E$ with $v\in E\subseteq V$, involves ${n-1\choose k-1}\leq (n-1)^{k-1}$ sets to be tested.
And we check this for the vertices in $V$ one after the other, until we either find one such $v$
or conclude none exists, so for at most $n$ vertices. So we test at most $n(n-1)^{k-1}$ sets.
Thus, by induction, we obtain the polynomial bound
$$t(n)\ \leq\ t(n-1)+n(n-1)^{k-1}\ \leq\
(n-1)^{k+1}+n(n-1)^{k-1}\ =\  (n^2-n+1)(n-1)^{k-1}\ \leq\ n^{k+1}\ .$$

As pointed out by one of the referees, the time bound above can be improved by
incorporating a suitable data structure, but as it does not affect the statement of the theorem,
and we wish to keep the article short and self contained, we do not elaborate on this.

\vskip.2cm{\em Proof of part 2.}
Now suppose $\H$ is orderable.
Assume that $n=km$ is a multiple of $k$ else there is no perfect matching.
We can then find an elimination order $v_1,\dots,v_n$ of $V$ in polynomial time by the
algorithm in the above proof of the first statement of the theorem. By definition each vertex $v_j$
for $j<k$ is both dominating and isolating and we choose to designate all of them as dominating.
Using this order we compute in polynomial time a sequence of integers $r_n,\dots,r_1$, initializing
$r_{n+1}:=0$ and setting $r_j:=r_{j+1}+k-1$ if $v_j$ is dominating and $r_j:=r_{j+1}-1$ if $v_j$ is
isolating for $j=n,\dots,1$.
Note that, letting $d_j$ and $i_j$ for $j=1,\dots,n$ be, respectively, the number of
dominating and isolating vertices with index at least $j$, we have that $r_j=(k-1)d_j-i_j$.

We claim $\H$ has a perfect matching if and only if all $r_i$ are nonnegative.
The intuition is that, when traversing the vertices backwards, $v_n,\dots,v_1$, each dominating
vertex $v_j$ can be used in a matching edge containing $k-1$ isolating vertices following it,
and therefore when encountering such a vertex we increase $r_j$ by $k-1$; on the other hand,
each isolating vertex $v_j$ takes a slot in such a matching edge and therefore when encountering
such a vertex we decrease $r_j$ by $1$.

But before proceeding with the proof of the claim, to illuminate it, following a suggestion
by one of the referees, we demonstrate how this claim leads to a process for finding a
perfect matching in an orderable hypergraph when all $r_i$ are nonnegative.
For $k=1$ the hypergraph $\H$ has a perfect matching if and only if
$\H=\{\{v_1\},\dots,\{v_n\}\}$ which holds if and only if all $r_i$ are nonnegative,
in which case $\M:=\H$ is the unique perfect matching.
So assume $k\geq 2$. Using the above elimination order $v_1,\dots,v_n$ and nonnegative
$r_1,\dots,r_n$, traverse the vertices backwards, $v_n,\dots,v_1$, and maintain two sets $D,I$,
of yet unmatched dominating and isolating vertices, respectively, and a set $\M\subseteq\H$ of a
partial matching constructed. The nonnegativity of the $r_i$ guarantees that each
$v_i$ in this backward order can be suitably matched.
Initialize $D:=I:=\emptyset$ and $\M:=\emptyset$. During the traversal proceed
as follows. If $v_i$ is dominating then update $D:=D\uplus\{v_i\}$. If $v_i$ is isolating then,
if $|I|<k-2$ then update $I:=I\uplus\{v_i\}$, whereas if $|I|=k-2$ then pick $v_r\in D$ with
largest index $r$, update $\M:=\M\uplus\{I\uplus\{v_i,v_r\}\}$ and update $D:=D\setminus\{v_r\}$
and $I:=\emptyset$. When the traversal is complete, add $k-|I|$ vertices from $D$ to $I$
including $v_r\in D$ with largest index $r$, and then add $I$ to $\M$. Finally, partition
the remaining vertices in $D$, if any, to sets of size $k$ and add all these sets to $\M$.
See Example \ref{process} below for a specific illustration of this process.

We now proceed with the proof of the claim that $\H$ has a perfect matching
if and only if all $r_i$ are nonnegative.
If $k=1$ then all $r_i$ are nonnegative if and only if all $v_i$ are
dominating which is equivalent to $\H=\{\{v_1\},\dots,\{v_n\}\}$,
which is indeed the case if and only if $\H$ has a perfect matching.
So assume $k\geq 2$. First we prove by induction on the number of isolating vertices that if
all $r_i$ are nonnegative then $\H$ has a perfect matching. Note that $v_n$ must be dominating else
$r_n=-1$ is negative. If there are $l\leq k-1$ isolating vertices then these vertices together
with $k-l$ dominating vertices including $v_n$ form an edge of $\H$. This edge, together with
$m-1$ more $k$-sets forming an arbitrary partition of the remaining $(m-1)k$ vertices,
which are all in $\H$ since all these vertices are dominating, forms a perfect matching.
Suppose then there are $l\geq k$ isolating vertices. Note that then $n\geq 2k$
since there are also at least $k$ dominating vertices $v_1,\dots,v_{k-1},v_n$.
Let $v_i$ have the $(k-1)$-largest index among the isolating vertices. Let $F$ be the $k$-set
consisting of $v_n$ and the $k-1$ isolating vertices of largest indices including $v_i$.
Note that $F\in\H$ since $v_n$ is dominating. Now let $V':=V\setminus F$ and
$\H':=\{E\in\H:E\subseteq V'\}$. It is clear that the order of $V'$ induced by the order
of $V$ is an elimination order for $H'$ with each vertex isolating or dominating as before.
Also $|V'|=n-k\geq k$ and the number of isolating vertices in $V'$ is smaller than that of $V$
by $k-1\geq 1$. Let $(r'_j:v_j\in V')$ be the sequence of integers computed using the elimination
order of $V'$. For each $v_j\in V$ let $d_j,i_j$, respectively, be the number of dominating and
isolating vertices in $V$ with indices at least $j$, and for each $v_j\in V'$ let $d'_j,i'_j$,
respectively, be the number of dominating and isolating vertices in $V'$ with indices at
least $j$. Note that for each $j$ we then have $r_j=(k-1)d_j-i_j$ and $r'_j=(k-1)d'_j-i'_j$.
Consider any $j$ with $v_j\in V'$. Note that $j\neq i,n$ since $v_i,v_n$ were removed from $V$.
If $j<i$ then $d'_j=d_j-1$ and $i'_j=i_j-(k-1)$ so $r'_j=r_j\geq 0$. If $j>i$ then all vertices
in $V'$ with indices at least $j$ must be dominating since all isolating vertices with indices
at least $j$ were removed from $V$, so $r'_j=(k-1)d'_j\geq k-1\geq 0$. Therefore all $r'_j$
are nonnegative and it follows by induction that $\H'$ has a perfect matching $\M'$.
Then $\M:=\M'\uplus\{F\}$ is a perfect matching of $\H$.

Second, suppose some $r_i$ is negative and let $D,I$, respectively, be the set of dominating and
isolating vertices in $V$ with indices at least $i$. Then $0>r_i=(k-1)|D|-|I|$ implies $|I|>(k-1)|D|$.
Suppose indirectly $\H$ has a perfect matching $\M$ and let $\M':=\{E\in\M: E\cap I\neq\emptyset\}$.
Since each $j\in I$ is isolating we must also have $E\cap D\neq\emptyset$ and hence
$|E\cap I|\leq k-1$ for all $E\in\M'$. This implies $|\M'|\leq|D|$ and
$(k-1)|\M'|\geq|I|$ which contradicts $|I|>(k-1)|D|$.
So if not all $r_i$ are nonnegative then $\H$ has no perfect matching. This completes the proof.
\epr

\be{process}
Here is an example of the process for constructing
a perfect matching in an orderable hypergraph described in the proof above.
Let $k=3$ and let $\H$ be the orderable $3$-hypergraph with $n=15$ vertices and elimination
order $v_1,\dots,v_{15}$ such that the vertices $v_1,v_2,v_3,v_4,v_6,v_{10},v_{13},v_{15}$ are
dominating and the vertices $v_5,v_7,v_8,v_9,v_{11},v_{12},v_{14}$ are isolating.
Going backwards we obtain the following sequence of the $r_i$ which are all nonnegative
and guarantee the existence of a perfect matching,
$$r_{15}=2,1,3,2,1,3,2,1,0,2,1,3,5,7,9=r_1\ .$$
Starting with $D:=I:=\emptyset$ and $\M:=\emptyset$ and traversing vertices backwards, we obtain:
\begin{eqnarray*}
v_{15} &:& D:=\{v_{15}\};\\
v_{14} &:& I:=\{v_{14}\};\\
v_{13} &:& D:=\{v_{13},v_{15}\};\\
v_{12} &:& D:=\{v_{13}\},\ I:=\emptyset,\ \M:=\{\{v_{12},v_{14},v_{15}\}\};\\
v_{11} &:& I:=\{v_{11}\};\\
v_{10} &:& D:=\{v_{10},v_{13}\};\\
v_9    &:& D:=\{v_{10}\},\ I:=\emptyset,\ \M:=\M\uplus\{\{v_9,v_{11},v_{13}\}\};\\
v_8    &:& I:=\{v_8\};\\
v_7    &:& D:=\emptyset,\ I:=\emptyset,\ \M:=\M\uplus\{\{v_7,v_8,v_{10}\}\};\\
v_6    &:& D:=\{v_6\};\\
v_5    &:& I:=\{v_5\};\\
v_4,v_3,v_2,v_1 &:& D:=\{v_1,v_2,v_3,v_4,v_6\},\ \M:=\M\uplus\{\{v_1,v_2,v_3\},\{v_4,v_5,v_6\}\}\ .
\end{eqnarray*}
So we end up with the perfect matching
$$\M\ =\ \{\{v_1,v_2,v_3\},\{v_4,v_5,v_6\},\{v_7,v_8,v_{10}\},
\{v_9,v_{11},v_{13}\},\{v_{12},v_{14},v_{15}\}\}\ .$$
\ee

\section{Separable hypergraphs}

We begin by comparing the classes of orderable hypergraphs and separable hypergraphs.
\bp{orderable_separable}
The following relations between orderable and separable hypergraphs hold:
\begin{enumerate}
\item
Every orderable $k$-hypergraph is separable;
\item
For $k=1,2$ every separable $k$-hypergraph
is orderable, but for each $k\geq 3$ there exist separable but not orderable $k$-hypergraphs.
\end{enumerate}
\ep

{\em Proof of part 1.}
Suppose $\H$ is orderable and let $v_1,\dots,v_n$ be an elimination order. Define a vertex
labeling by $a(v_i):=2^i$ if $v_i$ is dominating and $a(v_i):=-2^i$ if $v_i$ is isolating.
Consider any $k$-set $E\subseteq V$ and let $v_i\in E$ be the vertex with largest index in $E$,
which then implies that $v_i\in E\subseteq\{v_1,\dots,v_i\}$. If $v_i$ is dominating then $E\in\H$
and $a(E)\geq 2^i-\sum_{j<i}2^j>0$, whereas if $v_i$ is isolating then $E\notin\H$ and
$a(E)\leq -2^i+\sum_{j<i}2^j<0$.
Therefore we have that $\H=\{E\subseteq V:|E|=k,\ a(E)\geq 0\}$ is separable.

\vskip.2cm{\em Proof of part 2.}
For $k=1$, any $\H$ is orderable, since any order is an elimination order,
where we declare $v\in V$ isolating if $\{v\}\notin\H$ and $v\in V$ dominating if
$\{v\}\in\H$. For $k=2$, that is, graphs, we prove the claim by induction on $n:=|V|$.
For $n\leq 2$ any graph is orderable providing the induction base so assume $n\geq 3$.
Since $H$ is separable there exists a labeling $a:V\rightarrow\Z$ of vertices by integers such that
$$\H=\{\{u,v\}:u,v\in V,\ u\neq v,\ a(u)+a(v)\geq 0\}\ .$$
Now, either there is a $w\in V$ such that $a(w)\geq 0$ and $a(w)\geq|a(v)|$ for all $v\in V$,
or there is a $w\in V$ such that $a(w)<0$ and $|a(w)|>a(v)$ for all $v\in V$.
In either case, define $V':=V\setminus\{w\}$ and $\H':=\{E\in\H:E\subseteq V'\}$.
Clearly, for any two distinct $u,v\in V'$ we have that $\{u,v\}\in\H'$ if and only if
$a(u)+a(v)\geq 0$, so $\H'$ is separable, and $|V'|=n-1$, so by induction $\H'$ is orderable.
Let $v_1,\dots,v_{n-1}$ be an elimination order for $\H'$. We claim $v_1,\dots,v_{n-1},w$
is an elimination order for $\H$. Indeed, if $a(w)\geq 0$ then for all $v\in V'$ we have
$a(w)+a(v)\geq 0$ so we can declare $w$ dominating, whereas if $a(w)<0$ then for all $v\in V'$
we have $a(w)+a(v)<0$ so we can declare $w$ isolating. So the induction follows and we are done.

For $k\geq 3$ let $V=\{v_1,\dots,v_{k+1}\}$ and define a labeling $a:V\rightarrow\Z$ by setting
$$a(v_1):=0,\quad a(v_2):=\cdots:=a(v_k):=1,\quad a(v_{k+1}):=-(k-1)\ .$$
Then it is easy to check that the corresponding separable hypergraph is
$$\H=\{E\subseteq V:|E|=k,\ a(E)\geq 0\}\ =\
\{\{v_1,v_2,\dots,v_k\},\{v_2,v_3,\dots,v_{k+1}\}\}\ .$$
Now, suppose there exists an elimination order for $\H$ and let $v\in V$ be the last in that order.
If $v$ is isolating then it lies in no edge of $\H$. But this is impossible since every vertex
is in at least one of the two edges of $\H$. If $v$ is dominating then all $k$-sets $E$ with
$v\in E\subset V$ must be edges of $\H$. But there are $k\geq 3$ such sets whereas $\H$ has only
$2$ edges. So $v$ can be neither isolating nor dominating, hence there is no elimination
order and $\H$ is not orderable.
\epr

We remark, following a suggestion by one of the referees, that the construction in the
above proof can be extended to produce, for any $k\geq 3$ and any $n\geq k+1$, separable
but not orderable $k$-hypergraphs of order $n$, by simply adding $n-(k+1)$ more vertices
and labeling each by $-k$. The reader can verify that this works.

\vskip.2cm
We now prove our theorem about separable hypergraphs.

\vskip.2cm\noindent{\bf Theorem \ref{separable}}
{\em For any fixed $k$, the following two statements hold:
\begin{enumerate}
\item
It can be decided in polynomial time if a given $k$-hypergraph is separable;
\item
Deciding if a separable $k$-hypergraph has a perfect matching is
polynomial time doable for $k=1,2$, but is NP-complete for each $k\geq 3$, even
if $a:V\rightarrow\Z$ is encoded in unary.
\end{enumerate}}

\vskip.2cm{\em Proof of part 1.}
Given a $k$-hypergraph $\H$ on $V$, consider the following system
of linear inequalities in variables $a(v)$, $v\in V$,
$$\sum\{a(v):v\in E\}\ \geq\ 0\ \ \mbox{for each}\ \ E\in \H\ ,$$
$$\sum\{a(v):v\in E\}\ \leq\ -1\ \ \mbox{for each}\ \ E\subseteq V\ ,\ |E|=k\ ,\ E\notin \H\ .$$
Clearly $\H$ is separable if and only if this system has a rational solution. Since $k$ is fixed,
the number of inequalities is $O(n^k)$, polynomial in $n:=|V|$, and so the existence of a rational
solution of this system can be tested in polynomial time by linear programming, see e.g. \cite{Sch}.

\vskip.2cm{\em Proof of part 2.}
Now suppose $\H$ is separable. If $k=1,2$ then $\H$ is also orderable by Proposition
\ref{orderable_separable}, and so by Theorem \ref{orderable} the existence of a
perfect matching can be decided in polynomial time.
Suppose $k=3$. We reduce the NP-complete $3$-partition problem with input
$\H^==\{E\subseteq V:|E|=3,\ a(E)=0\}$ with $a$ encoded in unary to our problem over
$\H^\geq:=\{E\subseteq V:|E|=3,\ a(E)\geq 0\}$. We may and do assume $a(V)=0$ else $\H^=$
has no perfect matching. Indeed, if $E_1,\dots,E_m$ is a perfect matching of $\H^=$ then
$$a(V)\ =\ a\left(\uplus_{i=1}^m E_i\right)\ =\ \sum_{i=1}^m a(E_i)\ =\ 0\ .$$
We claim that $\H^=$ has a perfect matching if and only if $\H^\geq$ has.
Since $\H^=\subseteq\H^\geq$, clearly if $\H^=$ has a perfect matching then so does $\H^\geq$.
Conversely, suppose that $E_1,\dots,E_m$ is a perfect matching of $\H^\geq$. Then
$$\sum_{i=1}^m a(E_i)\ =\ a\left(\uplus_{i=1}^m E_i\right)\ =\ a(V)\ =\ 0\ ,$$
and $a(E_i)\geq 0$ for all $i$, so in fact $a(E_i)=0$ hence $E_i\in\H^=$ for all $i$.
Therefore $E_1,\dots,E_m$ is a perfect matching of $\H^=$ as well. So $3$-partition reduces to
perfect matching over separable $3$-hypergraphs.

Now suppose $k\geq 4$. We reduce the problem of deciding the existence of a perfect matching in a given
separable $3$-hypergraph $\H=\{E\subseteq V:|E|=3,\ a(E)\geq0\}$ where $a$ is encoded in unary,
just proved to be NP-complete, to the analog problem for separable $k$-hypergraphs.
We may and do assume $n:=|V|=3m$ is a multiple of $3$ and $\H\neq\emptyset$
else $\H$ has no perfect matching. Let $b:=1+\max\{a(E):E\subseteq V,|E|\leq k\}$
which can be computed in polynomial time as $k$ is fixed. Let $V':=V\uplus U$ where $U$
is a set of $(k-3)m$ new vertices. Define $a':V'\rightarrow\Z$ by $a'(u):=3b$ for $u\in U$
and $a'(v):=ka(v)-(k-3)b$ for $v\in V$. Let $\H':=\{E\subseteq V':|E|=k,\ a'(E)\geq0\}$.

We claim that $\H$ has a perfect matching if and only if $\H'$ has. Suppose first $E_1,\dots,E_m$ is a
perfect matching of $\H$. Let $U=\uplus_{i=1}^m F_i$ be an arbitrary partition of $U$ into $(k-3)$-sets
and let $E'_i:=E_i\uplus F_i$ for all $i$. Then $E'_1,\dots,E'_m$ is a perfect matching of $\H'$
since $V'=\uplus_{i=1}^m E'_i$, and for all $i$,
$$a'(E'_i)\ =\ \sum_{v\in E_i}(ka(v)-(k-3)b)+\sum_{u\in F_i}3b
\ =\ ka(E_i)-3(k-3)b+(k-3)3b\ =\ ka(E_i)\geq\ 0\ .$$
For the converse, consider first any $E'\in\H'$ and let $E:=E'\cap V$ and $F:=E'\cap U$. Then
\begin{eqnarray*}
0&\leq& a'(E')\ =\ \sum_{v\in E}(ka(v)-(k-3)b)+\sum_{u\in F}3b
\ =\ ka(E)-\sum_{v\in E}kb+\sum_{v\in E\uplus F}3b\\
& =& ka(E)-kb|E|+\sum_{v\in E'}3b\ \leq\ k(b-1)-kb|E|+k3b\ =\ kb(4-|E|)-k\ .
\end{eqnarray*}
This implies $1\leq b(4-|E|)$ and, since our assumption $\H\neq\emptyset$ implies $b\geq 1$,
we get $|E|\leq 3$.

Now suppose $E'_1,\dots,E'_m$ is a perfect matching of $\H'$. Let $E_i:=E_i'\cap V$ and $F_i:=E_i'\cap U$
for each $i$. Then $V=\uplus_{i=1}^mE_i$ hence $3m=|V|=\sum_{i=1}^m|E_i|$, and as just
proved, $|E_i|\leq 3$ for each $i$, so in fact $|E_i|=3$ and $|F_i|=k-3$ for each $i$.
Therefore, for each $i$, we find that $E_i\in\H$ since
$$0\ \leq\ a'(E'_i)\ =\ \sum_{v\in E_i}(ka(v)-(k-3)b)+\sum_{u\in F_i}3b
\ =\ ka(E_i)-3(k-3)b+(k-3)3b\ =\ ka(E_i)\ .$$
So $E_1,\dots,E_m$ is a perfect matching of $\H$.
Therefore the problem over $k$-hypergraphs is NP-complete for all $k\geq3$, completing the proof.
\epr

\section*{Acknowledgments}

The author was supported by a grant from the Israel Science Foundation and by the Dresner chair.
He thanks the referees for useful suggestions that improved the presentation.

\end{document}